\documentclass[12pt,a4paper,fleqn]{article}
\usepackage{a4wide,amsfonts,amsmath,latexsym,amssymb,euscript,graphicx,units,mathrsfs}

\usepackage{graphicx}
\usepackage{color}
\usepackage{amssymb}
\usepackage{amssymb}
\usepackage[T1]{fontenc}
\usepackage{latexsym}
\usepackage{xypic}
\usepackage{eufrak}
\usepackage{euscript}
\usepackage{amsfonts,amsmath}
\usepackage{verbatim}
\usepackage{fancyhdr}
\usepackage[english]{babel}
\usepackage{mathrsfs}
\usepackage{units}
\usepackage{cite}
\newtheorem{prop}{Proposition}[section]

\newtheorem{lemma}[prop]{Lemma}
\newtheorem{cor}[prop]{Corollary}

\newtheorem{theorem}[prop]{Theorem}
\newtheorem{thm}[prop]{Theorem}
\newtheorem{defi}[prop]{Definition}

\renewcommand{\geq}{\geqslant}
\def\leq{\leqslant}

\newcommand{\R}{\mathbb{R}}

\def\e{\varepsilon}
\font\tenbb=msbm10

\def\1{{\mathbf{1}}}

\font\tenbb=msbm10

\def\nN{\hbox{\tenbb N}}

\def\1{{\mathbf{1}}}
\def\0.5{{\frac{1}{2}}}

\newcommand{\fin}
{ \vspace{-0.6cm}
\begin{flushright}
\mbox{$\Box$}
\end{flushright}
\noindent }

\newcommand{\qed}{\nopagebreak\hspace*{\fill}
{\vrule width6pt height6ptdepth0pt}\par}

\newcounter{rea}
\begin{document}

\begin{center}
{\Large{\bf An invariance principle under the total variation distance}}
\normalsize
\\~\\ Ivan Nourdin\footnote{supported in part by the (french) ANR grant `Malliavin, Stein and Stochastic Equations with Irregular Coefficients'
[ANR-10-BLAN-0121].} (Universit\'e de Lorraine)\\
 Guillaume Poly (Universit\'e du Luxembourg)\\
\end{center}

{\small \noindent {{\bf Abstract}:
Let $X_1,X_2,\ldots$ be a sequence of i.i.d. \!\!random variables, with mean zero and variance one. Let $W_n=(X_1+\ldots+X_n)/\sqrt{n}$. An old and celebrated result of Prohorov \cite{prohorov} asserts that $W_n$
converges in  total variation to the standard Gaussian distribution if and only if $W_{n_0}$ has an absolutely continuous component for some $n_0$. In the present paper, we give yet another proof and extend Prohorov's theorem to a situation where, instead of $W_n$, we consider more generally a sequence of homogoneous polynomials in the $X_i$.
More precisely, we exhibit conditions for a recent invariance
principle proved by Mossel, O'Donnel and Oleszkiewicz \cite{MOO} to hold under the total variation distance.
There are many works about CLT under various metrics in the literature, but the present one
seems to be the first attempt to deal with homogeneous polynomials in the $X_i$ with degree {\it strictly} greater than one.
\\



\section{Introduction and main results}

Let $X_1, X_2, \ldots$ be independent copies of
a random variable with mean zero and variance one.
According to the central limit theorem,
the normalized sums
\begin{equation}\label{tcl}
W_n=\frac{X_1+\ldots+X_n}{\sqrt{n}}
\end{equation}
converge in distribution to the standard normal law $N\sim N(0,1)$.
In fact, using, e.g., the second Dini's theorem it is straightforward to
prove a much stronger result, namely that $W_n$ converges to $N$ in the Kolmogorov distance:
\begin{equation}\label{cvdkol}
\lim_{n\to\infty}\,\,d_{Kol}(W_n,N)=0,
\end{equation}
where
$d_{Kol}(U,V)=\sup_{x\in\R} \left| P(U\leq x)- P(V\leq x)\right|$.
\bigskip

In (\ref{cvdkol}), can we replace the Kolmogorov distance $d_{Kol}$ by the total variation distance $d_{TV}$, defined as
$d_{TV}(U,V)=\sup_{A\in\mathcal{B(\R)}} \left| P(U\in A)- P(V\in A)\right|$? In other words, do we also have
\begin{equation}\label{cvdtv}
\lim_{n\to\infty}\,\,d_{TV}(W_n,N)=0
\end{equation}
for $W_n$ defined by (\ref{tcl})?
The right answer is provided by an old and celebrated result of Prohorov \cite{prohorov}. To formulate it, first let us
introduce the Lebesgue decomposition (of the distribution) of a random variable. As is well known, each cumulative distribution function (cdf) $F$
may indeed be represented in the form:
\begin{equation}\label{lebesgue}
F(x)=u \int_{-\infty}^x g(y)dy+(1-u)G(x),\quad x\in\R,
\end{equation}
where $u\in[0,1]$, $g:\R\to [0,\infty)$ satisfies $\int_\R g(y)dy=1$ and $G$ is a singular cdf (corresponding to a distribution concentrated on a set of zero Lebesgue measure) with $G'(x)=0$ for almost all $x$. The real number $u\in[0,1]$ is uniquely determined by $F$; the density function $g$
is uniquely determined (up to a set of measure zero) if and only if $u>0$.

\begin{defi}
When $X$ is a random variable with cdf $F$, we say that $X$ is {\bf singular}
if $u=0$ in (\ref{lebesgue}).
If $u>0$, we say that $X$ has an {\bf absolutely continuous component} with density $g$.
\end{defi}

We can now state Prohorov's theorem \cite{prohorov}. A proof will be given in Section \ref{sec-pro}, to illustrate a possible use of our forthcoming results.

\begin{theorem}[Prohorov]\label{pro-thm}
The convergence (\ref{cvdtv}) takes place if and only if there exists $n_0\geq 1$ such that the random variable $W_{n_0}$ has an absolutely continuous component.
\end{theorem}
Prohorov's theorem has been the starting point of a fruitful line of research around the validity of the central limit theorem under various metrics and the estimation of their associated rates of convergence. Let us only give a small sample of references dealing with this rich and well studied topic. Convergence of densities in $L^\infty$ are studied by Gnedenko and Kolmogorov \cite{gne-kol}. On their side, Mamatov and Halikov \cite{Mamat-Vec} dealt with the multivariate CLT in total variation. Barron \cite{barron} studied the convergence in relative entropy, whereas Shimizu \cite{shimi} and Johnson and Barron \cite{joh-bar} studied the convergence in Fisher information. As far as rates of convergence are concerned, one can quote Mamatov and Sirazdinov \cite{Mamat} for the total variation distance and, very recently, Bobkov, Chistyakov and G\"{o}tze for bounds in entropy \cite{bobkov-B.E.E}, in Fisher information \cite{bobkov-F} and for Edgeworth-type expansions in the entropic central limit theorem \cite{bobkov-E}. Finally we mention \cite{Cac1, Cac2} for a variational approach of these issues with some variance bounds.

All these above-mentioned references have in common to `only' deal with {\it sums} of independent random variables. In the present paper, in contrast, we will consider highly {\it non-linear} functionals of independent random variables. It is arguably a much harder framework to work with, precisely because all the nice properties enjoyed by sums of independent variables are no longer valid in this context
(in particular, the use of characteristic functions is not appropriate).

Let us now turn to the details of the situation we are considering in the present article. Fix a degree of multilinearity $d\geq 1$
($d=1$ for linear, $d=2$ for quadratic, etc.) and, for any $n\geq 1$, consider a homogeneous polynomial  $Q_n:\R^{N_n}\to\R$ of the form
\begin{equation}\label{p1}
Q_n({\bf x})=\sum_{i_1,\ldots,i_d=1}^{N_n} a_n(i_1,\ldots,i_d)x_{i_1}\ldots x_{i_d},
\quad {\bf x}=(x_1,\ldots,x_{N_n})\in\R^{N_n}.
\end{equation}
In (\ref{p1}), it is implicitely supposed that $N_n\to\infty$ and also that $a_n(i_1,\ldots,i_d)$ are real numbers vanishing on diagonals and symmetric in the indices.
We further assume, in all what follows, that $Q_n$ is properly normalized:
\begin{equation}\label{p2}
d!\sum_{i_1,\ldots,i_d=1}^{N_n} a_n(i_1,\ldots,i_d)^2=1,\quad n\geq 1,
\end{equation}
and that all the terms of $Q_n$ are asymptotically negligible, meaning in our context that
\begin{equation}\label{nodominant}
\lim_{n\to\infty}\quad \max_{1\leq i_1\leq n}\sum_{i_2,\ldots,i_d=1}^{N_n} a_n(i_1,\ldots,i_d)^2\,\,= 0.
\end{equation}
As anticipated, property (\ref{nodominant}) will play a crucial role in the sequel.
It is also the key to obtain the following invariance principle, due to Mossel, O'Donnel and Oleszkiewicz \cite{MOO}.

\begin{thm}[Mossel, O'Donnel, Oleszkiewicz]\label{MOO}
Fix an integer $d\geq 1$, and let $Q_n$ be a sequence of homogeneous polynomials satisfying both (\ref{p1}),
(\ref{p2}) and
(\ref{nodominant}).
Let
${\bf X}=(X_1,X_2,\ldots)$  be
a sequence of independent
random variables with mean zero and variance one, belonging to $L^{2+\epsilon}(\Omega)$ for some $\epsilon>0$ (the same $\epsilon$ for each $X_i$).
Assume the same for
${\bf Y}=(Y_1,Y_2,\ldots)$. Then
\begin{equation}\label{cvdkol2}
\lim_{n\to\infty}\,\,d_{Kol}(Q_n({\bf X}),Q_n({\bf Y})) \,\,= 0.
\end{equation}
\end{thm}

Observe that one recovers (\ref{cvdkol}) by considering, in (\ref{cvdkol2}), $d=1$, $a_n(i)=\frac{1}{\sqrt{n}}$, $1\leq i\leq n$
(which satisfies (\ref{nodominant})) and $Y_1\sim N(0,1)$
(which leads to $Q_n({\bf Y})\sim N(0,1)$ for any $n$).\\

In the light of the aforementioned results, it seemed natural to us to ask under which assumption the convergence (\ref{cvdkol2}) may be strenghtened to the total variation distance as follows:
\begin{equation}\label{cvdtv2}
\lim_{n\to\infty}\,\,d_{TV}(Q_n({\bf X}),Q_n({\bf Y})) \,\,= 0.
\end{equation}

Before detailing our answer, let us first do
a quick digression.
As anticipated, a main aspect of our approach will consist in introducing
the following class of random variables.

\begin{defi}
For any $p\in]0,1]$ and $\alpha>0$, the class $\mathcal{C}(p,\alpha)$ is the set of
real random variables $X$ satisfying
\begin{equation}\label{law}
X\overset{\rm law}{=} \e(\alpha U+x_0)+(1-\e)V,
\end{equation}
where $x_0\in\R$ is a real number, and  $U\sim\mathcal{U}_{[-1,1]}$, $\varepsilon\sim \mathcal{B}(p)$ and $V$ (with no specified distribution) are three independent random variables.

\end{defi}

Since it will play a crucial role in the sequel, let us first
try to catch the meaning of (\ref{law}).
To this aim, we introduce
yet another class of random variables.
\begin{defi}
For any $c,\alpha>0$, the class $\mathcal{G}(c,\alpha)$ is the set of
real random variables $X$ having an absolutely continuous component and whose density $g$, see (\ref{lebesgue}), satisfies $g(x)\geq c$ for all $x\in[x_0-\alpha,x_0+\alpha]$ for some $x_0\in\R$.
\end{defi}

The following result  compares the two classes
$\mathcal{C}(p,\alpha)$
and $\mathcal{G}(c,\alpha)$.
Roughly speaking, it asserts that the class of random variables with an absolutely continous component (that is,
exactly the kind of random variables appearing in Prohorov's Theorem \ref{pro-thm}) coincides with
$\cup_{p\in]0,1],\alpha>0}\,\mathcal{C}_{p,\alpha}$.
Observe also that $\mathcal{G}(c,\alpha)$ is not empty if and only if $2c\alpha\leq 1$.

\begin{prop}\label{classC}
Fix $c,\alpha>0$ and $p\in]0,1]$.
One has $\mathcal{G}(c,\alpha)
\subset
\mathcal{C}(2c\alpha,\alpha)$.
Moreover, any random variable belonging to $\mathcal{C}(p,\alpha)$
has an absolutely continuous part.
\end{prop}

In Lemma \ref{aux} below, we state two further important properties of $\mathcal{C}(p,\alpha)$. Firstly, the sum of two independent random variables having an absolutely continuous component belong to $\cup_{c,\alpha>0}\,\mathcal{G}(c,\alpha)\subset\cup_{p\in]0,1],\alpha>0}\,\mathcal{C}(p,\alpha)$.
Secondly,  if $0<q\leq p\leq 1$ and $0<\beta\leq \alpha$, then
$\mathcal{C}(p,\alpha)\subset \mathcal{C}(q\beta/\alpha,\beta)$.\\

Now $\mathcal{C}(p,\alpha)$ has been introduced and is arguably well-understood, let us  give a name to the set of sequences of independent and normalized random variables we will deal with throughout the sequel.

\begin{defi}
Let $\alpha>0$, $p\in]0,1]$ and $\epsilon>0$. A sequence ${\bf X}=(X_1,X_2,\ldots)$
of random variables belongs to $\mathcal{D}(\alpha,p,2+\epsilon)$
if the $X_i$ are independent, satisfy $\sup_i E|X_i|^{2+\epsilon}<\infty$ and if, for each $i$,  $E[X_i]=0$, $E[X_i^2]=1$ and $X_i\in\mathcal{C}(p,\alpha)$.
\end{defi}

We are now in a position to state the main result of the present paper.
\begin{thm}\label{main}
Fix an integer $d\geq 1$, and let $Q_n$ be a sequence of homogeneous polynomials satisfying both (\ref{p1}), (\ref{p2} and
(\ref{nodominant}).
Let ${\bf X}$ and ${\bf Y}$ belong to $\mathcal{D}(\alpha,p,2+\epsilon)$ for some $\epsilon,\alpha>0$ and $p\in]0,1]$.
Then (\ref{cvdtv2}) holds true.
\end{thm}

A noticeable corollary of Theorem \ref{main} is a new proof of Prohorov's Theorem \ref{pro-thm}, see Section \ref{sec-pro}.
Another one is the following result.
\begin{cor}\label{cor1}
Fix an integer $d\geq 2$, and let $Q_n$ be a sequence of homogeneous polynomials satisfying (\ref{p1})-(\ref{p2})-(\ref{nodominant}).
Let ${\bf X}$ belong to  $\mathcal{D}(\alpha,p,2+\epsilon)$ for some $\epsilon,\alpha>0$ and $p\in]0,1]$.
 If $Q_n({\bf X})$ converges in law to $W$,
then $W$ has a density and $Q_n({\bf X})$ converges to $W$ in total variation.
\end{cor}
The statement of Corollary \ref{cor1} would be clearly wrong without assuming (\ref{nodominant}). Consider, e.g., $Q_n({\bf x})=x_1$, $n\geq 1$ with $X_1$ singular. Another interesting consequence
of Theorem \ref{main} is provided by the next theorem.

\begin{thm}\label{CLT}
Let $\{a_n(i_1,\ldots,i_d)\}_{1\leq i_1,\ldots,i_d\leq N_n}$ be an array of real numbers vanishing on diagonals, symmetric in the indices
and satisfying (\ref{p2}).
(We do \underline{not} suppose (\ref{nodominant}).)
Let ${\bf G}=(G_1,G_2,\ldots)$ be a sequence of independent $N(0,1)$ random variables, and
let $N\sim N(0,1)$.
Then, the following four assertions are equivalent as $n\to\infty$.
\begin{itemize}
\item[(a)] $\displaystyle{
\sum_{i_1,i_2,\cdots,i_d=1}^{N_n} a_n (i_1,i_2,\cdots,i_d) G_{i_1} G_{i_2} \cdots G_{i_d}\overset{\rm law}{\to} N}$.
\item[(b)] $\displaystyle{
d_{TV}\left(N,\sum_{i_1,i_2,\cdots,i_d=1}^{N_n} a_n (i_1,i_2,\cdots,i_d) G_{i_1} G_{i_2} \cdots G_{i_d}\right)}\to 0$.
\item[(c)] For all ${\bf X}=(X_1,X_2,\ldots)$ belonging to  $\mathcal{D}(\alpha,p,2+\epsilon)$
for some $\epsilon,\alpha>0$ and $p\in]0,1]$, we have
\begin{equation*}
\sum_{i_1,i_2,\cdots,i_d=1}^{N_n} a_n (i_1,i_2,\cdots,i_d) X_{i_1} X_{i_2} \cdots X_{i_d}\overset{\rm law}{\to} N.
\end{equation*}
\item[(d)] For all ${\bf X}=(X_1,X_2,\ldots)$ belonging to  $\mathcal{D}(\alpha,p,2+\epsilon)$
for some $\epsilon,\alpha>0$ and $p\in]0,1]$, we have
\begin{equation*}
d_{TV}\left(N,\sum_{i_1,i_2,\cdots,i_d=1}^{N_n} a_n (i_1,i_2,\cdots,i_d) X_{i_1} X_{i_2} \cdots X_{i_d}\right)\to 0.
\end{equation*}
\end{itemize}
\end{thm}

The rest of our paper is organised as follows. In Section 2 we prove all the results that are stated in this Introduction, except Theorem \ref{main}; in particular, we give a new proof of Prohorov's Theorem \ref{pro-thm}  in Section \ref{sec-pro}.
Finally, the proof of our main result, namely Theorem \ref{main}, is  provided in Section 3.

\section{Proofs of all stated results except Theorem \ref{main}}

\subsection{Some useful lemmas}

The following lemma will be used several times in the sequel.
\begin{lemma}\label{riesz}
Fix $q\in[0,1]$, and
let $Y,Z$ be two random variables
satisfying $E[f(Y)]\geq qE[f(Z)]$ for all positive bounded function $f$. Then there exists two independent random variable $W$ and $\zeta\sim\mathcal{B}(q)$, independent from $Z$, such that
\begin{equation}\label{id-law}
Y\overset{\rm law}{=} \zeta Z +(1-\zeta)W.
\end{equation}
\end{lemma}
{\it Proof of Lemma \ref{riesz}}.
Our assumption ensures that the linear form $f\mapsto
E[f(Y)]- qE[f(Z)]$ is positive.
From the Riesz representation theorem, one deduces the existence of a positive finite Radon measure $\nu$
such that
\begin{equation}\label{radon}
E[f(Y)]= qE[f(Z)]+\int_\R f(x)d\nu(x).
\end{equation}
Choosing $f\equiv 1$ in (\ref{radon}) gives
$\nu(\R)=1-q$.
If $\nu(\R)=0$ then $q=1$ and the proof of (\ref{id-law}) is established.
Otherwise, $\nu(\R)>0$ and one can consider $W\sim \frac{1}{\nu(\R)}d\nu(x)$, implying in turn (\ref{id-law}).\fin

In the following lemma, we gather useful properties of the
classes $\mathcal{C}(p,\alpha)$ and $\mathcal{G}(c,\alpha)$.
\begin{lemma}\label{aux}
The following properties take place.
\begin{enumerate}
\item  If $0<q\leq p\leq 1$ and if $0<\beta\leq \alpha$, then
$\mathcal{C}(p,\alpha)\subset \mathcal{C}(q\beta/\alpha,\beta)$. In particular, $\mathcal{C}(p,\alpha)\subset \mathcal{C}(q,\alpha)$.
\item If $X$ and $Y$ both have an absolutely continuous component and if $X$ is independent from $Y$, then there exists $c,\alpha>0$ such that $X+Y\in\mathcal{G}(c,\alpha)$.
\item If $X$ belongs to $\mathcal{C}(p,\alpha)$
with $\alpha>0$ and $p\in]0,1]$ and if $Y$ is any random variable independent from $X$, then $X+Y$ belongs to
$\mathcal{C}(q,\beta)$
for some $\beta>0$ and $q\in]0,1]$.
\item If $a\neq 0$ and $b$ are two real numbers
and if $X$ belongs to $\mathcal{C}(p,\alpha)$
with $\alpha>0$ and $p\in]0,1]$, then
$aX+b\in\mathcal{C}(p,|a|\alpha)$.
\end{enumerate}
\end{lemma}
{\it Proof}. 1.
Fix  $0<q\leq p\leq 1$ and $0<\beta\leq \alpha$, and consider
$X\in\mathcal{C}(p,\alpha)$. According to (\ref{law}), we have,
for any positive $f$,
\begin{eqnarray*}
E[f(X)]&=&p\int_\R f(x)\frac{1}{2\alpha}{\bf 1}_{[x_0-\alpha,x_0+\alpha]}(x)dx+(1-p)E[f(V)]\\
&\geq&\frac{q\beta}{\alpha}\int_\R f(x)\frac{1}{2\beta}{\bf 1}_{[x_0-\beta,x_0+\beta]}(x)dx.
\end{eqnarray*}
The conclusion follows from Lemma \ref{riesz}.\\
2. Consider the decomposition (\ref{lebesgue}) of the cdf $F$ of $X$. This settles $u\in]0,1]$ and $g:\R_+\to\R$ in a unique way.
Settle similarly $v\in]0,1]$ and $h:\R_+\to\R$ for $Y$.
For any Borel set $A$, one has $P(X\in A)\geq u\int_A g(x) dx$ and the same for $Y$. We deduce
\[
P(X+Y\in A)\geq  uv\int_A (g\star h) (x) dx,
\]
with $\star$ denoting the usual convolution. Besides, $g\star h=\lim_{M\to\infty} g \star \inf(h,M)$ and the limit is increasing by positivity of $g$. Finally, we note that, since $g\in L^1$ and $\inf(h,M)\in L^\infty$, the convolution $g\star \inf(h,M)$ is continuous. Let $x_0\in\R$ and $M>0$ be such that $(g\star \inf(h,M)) (x_0)>0$.
(Such a pair $(x_0,M)$ necessarily exists, otherwise we would have $g\star h\equiv 0$ by taking the large $M$ limit.) By continuity, there exists $c>0$ and $\alpha>0$ such that, for any $x\in]x_0-\alpha,x_0+\alpha[$, $(g\star h) (x)\geq (g \star \inf(h,M))(x) \ge c$. That is, $X+Y$ belongs to $\mathcal{G}(c,\alpha)$.\\
3.  We have $X\overset{\rm law}{=}\e(\alpha U+x_0)+(1-\e)V$, with
$x_0\in\R$ a real number, and  $U\sim\mathcal{U}_{[-1,1]}$, $\varepsilon\sim \mathcal{B}(p)$ and $V$ (with no specified distribution)  three independent random variables.
On the other hand, one can write
$Y\overset{\rm law}{=}\e Y+(1-\e)Z$, with $Z$ having the same law than $Y$ and independent from $Y,U,\e,V$. Thus,
\[
X+Y\overset{\rm law}{=}\e(\alpha U+x_0+Y)+(1-\e)(V+Z).
\]
The random variable $\alpha U+x_0+Y$ has a density
$g$ given by
\[
g(v)=\int_\R \frac{1}{2\alpha}{\bf 1}_{[x_0-\alpha,x_0+\alpha]}(v-y)dP_Y(y)=
\frac{1}{2\alpha}\,P(Y\in[v-x_0-\alpha,v-x_0+\alpha]).
\]
As a matter of fact, $g$ is a regulated function, since it is the difference of two increasing functions. In particular,
the set $\mathcal{E}$ of its discontinuous points is countable. As a consequence, ${\rm Leb}(\mathcal{E})=0$, implying in turn
$1=\int_\R g(v)dv = \int_{\R\setminus \mathcal{E}}g(v)dv$,
so that there exists $x_1\not\in\mathcal{E}$ satisfying
$g(x_1)>0$. Since $g$ is continuous at $x_1$, there exists $r>0$ such that $g(v)\geq \frac12g(x_1)$ for all
$v\in[x_1-r,x_1+r]$. By Lemma \ref{riesz}, it comes that
\[
\alpha U+x_0+Y\overset{\rm law}{=}\eta(rU+x_1)+(1-\eta)T,
\]
where $\eta\sim \mathcal{B}(p')$ for some $p'\in]0,1]$,
$U\sim \mathcal{U}_{[-1,1]}$ and $T$ are independent. Hence
\[
X+Y\overset{\rm law}{=}\e\eta(rU+x_1)+\e(1-\eta)T+(1-\e)(V+Z).
\]
As a result, for any bounded positive function,
\begin{eqnarray*}
E[f(X+Y)]&=&pp'E[f(rU+x_1)]+p(1-p')E[f(T)]+(1-p)E[f(V+Z)]\\
&\geq&pp'E[f(rU+x_1)].
\end{eqnarray*}
Finally, one deduces that $X+Y$ belongs to $\mathcal{C}(pp',r)$ by Lemma \ref{riesz}.\\
4. Obvious.
\fin

\subsection{Proof of Proposition \ref{classC}} \label{sec-prop}

Let $X$ be an element of $\mathcal{G}(c,\alpha)$.
Let $F$ denote its cdf, and consider $g$, $G$ and $u$ as
in (\ref{lebesgue}).
Let $Y\sim g(x)dx$, $Z\sim dG(x)$ and $\eta\sim\mathcal{B}(u)$ be three independent
random variables.
Then $X\overset{\rm law}{=}\eta Y+(1-\eta)Z$.
Using the assumption made on $g$, one obtains,
for any positive function $f$,
\[
E[f(Y)]\geq 2c\alpha \,E\left[f\left(\alpha U+x_0\right)\right].
\]
Observe that $0<2c\alpha\leq 1$ necessarily.
We then deduce that $X\in\mathcal{C}(2c\alpha,\alpha)$ from Lemma \ref{riesz}.

Consider now a random variable $X$ belonging to $\mathcal{C}(p,\alpha)$. We have, for any positive function $f$ and according to the decomposition (\ref{law}),
\begin{equation}\label{gus}
E[f(X)]=pE[f(\alpha U+x_0)]+(1-p)E[f(V)].
\end{equation}
Let us consider the Lebesgue decomposition $(u,g,G)$ of  $V$,
see (\ref{lebesgue}):
\[
E[f(V)]=u\int_\R f(x)g(x)dx + (1-u)\int_\R f(x)dG(x).
\]
Plugging into (\ref{gus}) yields
\[
E[f(X)]=\int_\R f(x)\left\{\frac{p}{2\alpha}{\bf 1}_{[x_0-\alpha,x_0+\alpha]}(x)
+(1-p)u\,g(x)
\right\}dx+(1-p)(1-u)\int_\R f(x)dG(x),
\]
from which we deduce that $X$ has an absolutely continuous component.
\fin

\subsection{Proof of Corollary \ref{cor1}}

Let the assumption of Corollary \ref{cor1} prevail
and consider a sequence ${\bf G}=(G_1,G_2,\ldots)$ composed of independent copies of a standard Gaussian random variable.
By the Mossel, O'Donnel, Oleszkiewicz's invariance principle (\ref{cvdkol2}), one has that $Q_n({\bf G})$ converges in law to $W$.
It implies, by \cite[Theorem 3.1]{NP} (see also \cite[Lemma 2.4]{NP}), that $W$ has a density and $d_{TV}(Q_n({\bf G}),W)\to 0$ as $n\to\infty$.
But $d_{TV}(Q_n({\bf X}),Q_n({\bf G}))\to 0$ by our Theorem \ref{main}, so
\[
d_{TV}(Q_n({\bf X}),W)\leq
d_{TV}(Q_n({\bf X}),Q_n({\bf G}))
+d_{TV}(Q_n({\bf G}),W)\to 0\quad\mbox{as $n\to\infty$}.
\]
This concludes the proof of Corollary \ref{cor1}.\fin

\subsection{Proof of Theorem \ref{pro-thm}}\label{sec-pro}

We want to use Corollary \ref{cor1} but the problem is that ${\bf X}$ is not assumed to belong to $\mathcal{D}$ in the statement of Theorem \ref{pro-thm}.
To overcome this issue we shall need Lemma \ref{aux}.

Let the assumptions and notation of Theorem \ref{pro-thm}
prevail.

First, when $W_{n}$ is singular, then there exists
a Borel set $A$ such that $P(N\in A)=0$ and $P(W_n\in A)=1$; in particular, $d_{TV}(W_n,N)=1$.
Hence, if $W_{n}$ is singular
for all $n$, then (\ref{cvdtv}) cannot hold.

Now, assume that $W_{n_0}$ has an absolutely continuous component for some $n_0$.
To prove (\ref{cvdtv}) is obviously equivalent to prove that,
for all $k\in\{0,\ldots,2n_0-1\}$,
\begin{equation}\label{cvdtvbis}
\lim_{n\to\infty}d_{TV}(W_{2n_0n+k}, N)=0.
\end{equation}
So, fix $k\in\{0,\ldots,2n_0-1\}$ and
consider a sequence $(Y_2,Y_3,\ldots)$ of independent copies of $W_{2n_0}$.
By Lemma \ref{aux} (points 2 and 4), observe that each $Y_i$
belongs to $\mathcal{C}(p,\alpha)$, for some $p\in]0,1]$ and
$\alpha>0$ (the same $p$ and the same $\alpha$ for all $i\geq 2$); also, we have $E[Y_i]=0$ and $E[Y_i^2]=1$.
On the other hand, let $Y_1$ be independent of $Y_2,Y_3,\ldots$ and have the same law than $W_{2n_0+k}$.
By Lemma \ref{aux} (points 3 and 4),
$Y_1$
belongs to $\mathcal{C}(q,\beta)$ for some $q\in]0,1]$ and
$\beta>0$; also, we have $E[Y_1]=0$ and $E[Y_1^2]=1$.
In fact, thanks to Lemma \ref{aux} (point 1), one may and will choose the same $p$ and the same $\alpha$ for each $Y_i$, without making a difference between $i=1$ and $i\geq 2$.

Bearing all the previous notation in mind, we can write
\[
W_n \overset{\rm law}{=}
\sqrt{\frac{2n_0+k}{n}}\,Y_1+
\sqrt{\frac{2n_0}{n}}\,\sum_{k=2}^{n}Y_k.
\]
The convergence (\ref{cvdtvbis}) is now a direct
consequence of Theorem \ref{main} applied to $d=1$, a sequence ${\bf X}={\bf G}=(G_1,G_2,\ldots)$ of independent $N(0,1)$ variables, $a_n(1)=\sqrt{\frac{2n_0+k}{n}}$ and $a_n(i)=\sqrt{\frac{2n_0}{n}}$, $i=2,\ldots,n$.
Note that
\[
d_{FM}(W_n,N)=d_{FM}(Q_n({\bf Y}),Q_n({\bf X}))\to 0\quad \mbox{as $n\to\infty$}
\]
by the usual CLT, so that it is not necessary to rely on Theorem \ref{MOO} to conclude the proof of Theorem \ref{main} (see Step 7 of Section \ref{sec-main}) and so to assume the existence of an absolute $q$th moment for $X_1$ with $q$ {\it strictly} greater than 2.\fin

\subsection{Proof of Theorem \ref{CLT}}

Implication $(a)\Rightarrow(b)$ is a reformulation of \cite[Corollary 5.2.8]{NouPecBook}.
Implications $(b)\Rightarrow(a)$ and $(d)\Rightarrow (c)$ are obvious. Implication $(c)\Rightarrow(a)$ is
because ${\bf G}$ belongs to $\mathcal{D}(\alpha,p,2+\e)$ for some $\alpha,\e>0$ and $p\in]0,1[$ (use Proposition \ref{classC}).
Finally, implication $(a)\Rightarrow (d)$ is a consequence of the following two facts.
Firstly, if $(a)$ takes place then, by a usual hypercontractivity argument,
the sequence $F_n:=\sum_{i_1,i_2,\cdots,i_d=1}^{N_n} a_n (i_1,i_2,\cdots,i_d) G_{i_1} G_{i_2} \cdots G_{i_d}$ (which is normalized so that $E[F_n^2]=1$, see indeed (\ref{p2})) satisfies $E[F_n^4]\to 3$.
Secondly, one has, according to \cite[(11.4.7) and (11.4.8) pp. 192-193]{NouPecBook}:
\[
\max_{1\leq i_1\leq n}\sum_{i_2,\ldots,i_d=1}^{N_n} a_n(i_1,\ldots,i_d)^2\leq \frac{1}{dd!}
\sqrt{E[F_n^4]-3}.
\]
These two facts together imply that, if $(a)$ holds, then (\ref{nodominant}) is automatically satisfied. Thus, Theorem \ref{main} implies the validity of assertion $(d)$.\fin

\section{Proof of Theorem \ref{main}}\label{sec-main}
Let the assumptions and notation of Theorem \ref{main}.
Without loss of generality, for simplicity we assume
that $N_n=n$.

The proof is divided into several steps.\\

{\bf Step 1}.
In the definition of $Q_n({\bf X})$ one may and will replace each $X_i$ by $\e_i(\alpha U_i+x_i)+(1-\e_i)V_i$,
where ${\bf e}=(\e_1,\e_2,\ldots)$ is a sequence of independent Bernoulli random variables ($\e_i\sim \mathcal{B}(p_i)$), ${\bf U}=(U_1,U_2,\ldots)$ is a sequence of independent $[-1,1]$-uniformly distributed random variables and  ${\bf V}=(V_1,V_2,\ldots)$ is a sequence of independent random variables; moreover, ${\bf e}$, ${\bf U}$ and ${\bf V}$ are independent. That is,
\[
Q_n({\bf X})=\sum_{i_1,\ldots,i_d=1}^n
a_n(i_1,\ldots,i_d)\{\e_{i_1}(\alpha_{i_1} U_{i_1}+x_{i_1})+(1-\e_{i_1})V_{i_1}\}\ldots
\{\e_{i_d}(\alpha_{i_d} U_{i_d}+x_{i_d})+(1-\e_{i_d})V_{i_d}\}.
\]

Now, let us expand everything, and then rewrite
$Q_n({\bf X})$ as a polynomial in the $U_i$. We obtain
\[
Q_n({\bf X})=A_n+B_n+C_n,
\]
where
\begin{eqnarray*}
A_n&=&\sum_{i_1,\ldots,i_d=1}^n
a_n(i_1,\ldots,i_d)\alpha_{i_1}\ldots\alpha_{i_d}
\e_{i_1}\ldots\e_{i_d}
U_{i_1}\ldots U_{i_d}\\
B_n&=&Q_n({\bf X})-A_n-C_n\\
C_n&=&\sum_{i_1,\ldots,i_d=1}^n ,
a_n(i_1,\ldots,i_d)
\{\e_{i_1}x_{i_1}+(1-\e_{i_1})V_{i_1}\}\ldots
\{\e_{i_d}x_{i_d}+(1-\e_{i_d})V_{i_d}\}
\end{eqnarray*}
satisfying
\[
E[A_nB_n|{\bf e},{\bf V}]=
E[A_nC_n|{\bf e},{\bf V}]=
E[B_nC_n|{\bf e},{\bf V}]=
E[A_n|{\bf e},{\bf V}]=
E[B_n|{\bf e},{\bf V}]=0.
\]
As a result,
\begin{eqnarray*}
{\rm Var}[Q_n({\bf X})|{\bf e},{\bf V}]&=&
E[A_n^2|{\bf e},{\bf V}]+E[B_n^2|{\bf e},{\bf V}]
+{\rm Var}[C_n|{\bf e},{\bf V}]\\
&\geq&E[A_n^2|{\bf e},{\bf V}]\\
&\geq&\alpha^{2d}3^{-d}d!
\sum_{i_1,\ldots,i_d=1}^n
a_n(i_1,\ldots,i_d)^2
\e_{i_1}\ldots\e_{i_d}.
\end{eqnarray*}
To go one step further, let us decompose $\e_i$ into
$(\e_i-p)+p$ and use (\ref{p2}), so to obtain
\begin{eqnarray*}
&&{\rm Var}[Q_n({\bf X})|{\bf e},{\bf V}]\\
&\geq&
d!\left(\frac{\alpha^{2}p}{3}\right)^d
+
\alpha^{2d}3^{-d}d!
\sum_{k=1}^d \frac{\binom{d}{k}}{p^{k}}
\sum_{i_1,\ldots,i_{d}=1}^n
a_n(i_1,\ldots,i_d)^2
(\e_{i_1}-p)\ldots(\e_{i_k}-p).
\end{eqnarray*}
Using the assumptions  (\ref{p2}) and (\ref{nodominant}),
we can write,
for any fixed $k\in\{1,\ldots,d\}$,
\begin{eqnarray*}
&&E\left[\left(
\sum_{i_1,\ldots,i_{d}=1}^n
a_n(i_1,\ldots,i_d)^2
(\e_{i_1}-p)\ldots(\e_{i_k}-p)
\right)^2\right]\\
&=&p^k(1-p)^k
\sum_{i_1,\ldots,i_{k}=1}^n\left(
\sum_{i_{k+1},\ldots,i_d=1}^n a_n(i_1,\ldots,i_d)^2
\right)^2\\
&\leq&\sum_{i_1=1}^n
\left(
\sum_{i_{2},\ldots,i_d=1}^n a_n(i_1,\ldots,i_d)^2
\right)^2\\
&\leq&\max_{1\leq j_1\leq n}
\sum_{j_{2},\ldots,j_d=1}^n a_n(j_1,\ldots,j_d)^2
\times
\sum_{i_{1},\ldots,i_d=1}^n a_n(i_1,\ldots,i_d)^2\\
&=&\frac{1}{d!}\,\,\,\max_{1\leq j_1\leq n}
\sum_{j_{2},\ldots,j_d=1}^n a_n(j_1,\ldots,j_d)^2\,\,\to 0
\quad\mbox{as $n\to\infty$}.
\end{eqnarray*}
We deduce that, in probability,
\begin{equation}\label{cl-step1}
\liminf_{n\to\infty} {\rm Var}[Q_n({\bf X})|{\bf e},{\bf V}]\geq
d!\left(\frac{\alpha^{2}p}{3}\right)^d.
\end{equation}

{\bf Convention}. From now on, and since all the quantities we are dealing with are measurable with respect to ${\bf e}$, ${\bf U}$ and ${\bf V}$,
we shall write $E_{\bf U}$ (resp. $E_{\bf e,V}$) to indicate the mathematical expectation with respect to ${\bf U}$ (resp. ${\bf e}$ and ${\bf V}$)
Note that $E_{\bf U}$ coincides with the conditional expectation $E[\cdot|{\bf e},{\bf V}]$.\\

{\bf Step 2}. Set $p_{\alpha}(x)=\frac{1}{\alpha\sqrt{2\pi}}e^{-\frac{x^2}{2\alpha^2}}$, $x\in\R$, $0<\alpha\leq 1$, and
let $\phi\in C^\infty_c$ be bounded by 1.
It is immediately checked that
\begin{equation}\label{fc1}
\|\phi\star p_\alpha\|_\infty\leq 1\leq\frac1\alpha
\quad\mbox{and}\quad
\|(\phi\star p_\alpha)'\|_\infty
\leq\frac1\alpha.
\end{equation}
We can write
\begin{eqnarray*}
&&\left|
E[\phi(Q_n({\bf X})]-E[\phi(Q_n({\bf Y})]
\right|\\
&\leq&
\left|E[(\phi-\phi\star p_\alpha)(Q_n({\bf X}))]
\right|
+
\left|E[(\phi-\phi\star p_\alpha)(Q_n({\bf Y}))]
\right|
+\frac{1}{\alpha}d_{FM}(Q_n({\bf X}),Q_n({\bf Y})),
\end{eqnarray*}
where $d_{FM}$ stands for the Fortet-Mourier distance, which is known to metrize the convergence in law.
Let us concentrate on the first two terms. We have, e.g., for the first term:
\begin{eqnarray*}
&&
\left|E[(\phi-\phi\star p_\alpha)(Q_n({\bf X}))]
\right|\\
&\leq&
\left|E[(\phi-\phi\star p_\alpha)(Q_n({\bf X}))
{\bf 1}_{\{
{\rm Var}[Q_n({\bf X})|{\bf e},{\bf V}]<
\frac{d!}{2}\left(\frac{\alpha^{2}p}{3}\right)^d
\}}
]
\right|\\
&&+
\left|E[(\phi-\phi\star p_\alpha)(Q_n({\bf X}))
{\bf 1}_{\{
E[Q_n({\bf X})^2|{\bf e},{\bf V}]>M
\}}
]
\right|\\
&&+\left|E[(\phi-\phi\star p_\alpha)(Q_n({\bf X}))
{\bf 1}_{\{
{\rm Var}[Q_n({\bf X})|{\bf e},{\bf V}]\geq
\frac{d!}{2}\left(\frac{\alpha^{2}p}{3}\right)^d,\,\,
E[Q_n({\bf X})^2|{\bf e},{\bf V}]\leq M
\}}
]
\right|\\
&\leq&2P\left(
{\rm Var}[Q_n({\bf X})|{\bf e},{\bf V}]<
\frac{d!}{2}\left(\frac{\alpha^{2}p}{3}\right)^d
\right)
+2P\left(E[Q_n({\bf X})^2|{\bf e},{\bf V}]>M\right)\\
&&+\left|E[(\phi-\phi\star p_\alpha)(Q_n({\bf X}))
{\bf 1}_{\{
{\rm Var}[Q_n({\bf X})|{\bf e},{\bf V}]\geq
\frac{d!}{2}\left(\frac{\alpha^{2}p}{3}\right)^d,\,\,
E[Q_n({\bf X})^2|{\bf e},{\bf V}]\leq M
\}}
]
\right|.
\end{eqnarray*}
We have, using the Markov inequality,
\[
P\left(E[Q_n({\bf X})^2|{\bf e},{\bf V}]>M\right)\leq
\frac{1}M\,E[E[Q_n({\bf X})^2|{\bf e},{\bf V}]]=\frac1M.
\]
On the other hand,
\begin{eqnarray*}
&&\left|E\left[(\phi-\phi\star p_\alpha)(Q_n({\bf X}))
{\bf 1}_{\{
{\rm Var}[Q_n({\bf X})|{\bf e},{\bf V}]\geq
\frac{d!}{2}\left(\frac{\alpha^{2}p}{3}\right)^d,\,\,
E[Q_n({\bf X})^2|{\bf e},{\bf V}]\leq M
\}}
\right]
\right|\\
&\leq&
E_{\bf e,V}\left[\big|E_{\bf U}[(\phi-\phi\star p_\alpha)(Q_n({\bf X}))
]
\big|
{\bf 1}_{\{
{\rm Var}[Q_n({\bf X})|{\bf e},{\bf V}]\geq
\frac{d!}{2}\left(\frac{\alpha^{2}p}{3}\right)^d,\,\,
E[Q_n({\bf X})^2|{\bf e},{\bf V}]\leq M
\}}
\right].
\end{eqnarray*}

{\bf Step 3}. In this step, we shall introduce the framework we are going to use for the rest of the proof.
We refer the reader to \cite{Ledoux-bakry-gentil} for the details and missing proofs.
Fix an integer $m$ and let $\mu$ denote the distribution of the random vector $(X_1,\ldots,X_{m})$,
with $X_1,\ldots,X_m$ independent copies of $U\sim \mathcal{U}_{[-1,1]}$,
There exists a reversible Markov process on $\R^m$, with semigroup $P_t$, equilibrium measure $\mu$ and generator $\mathcal{L}$
given by
\begin{equation}\label{Jac}
\mathcal{L}f(x)=\sum_{i=1}^{m} \Big{(}(1-x_{i}^2) \partial_{ii}f-2 x_i\,\partial_i f\Big{)},\quad x\in\R^m.
\end{equation}
The operator $\mathcal{L}$ is selfadjoint and
negative semidefinite. We define the carr\'e du champ operator $\Gamma$ as
\begin{equation}\label{gamma3}
\Gamma(f,g)(x) =
 \frac12\big(\mathcal{L}(fg)(x)-f(x)\mathcal{L}g(x)-g(x)\mathcal{L}f(x)\big)
=
\sum_{i=1}^{m} (1-x_i^2) \partial_i f(x)\partial_i g(x).
\end{equation}
When $f=g$ we simply write $\Gamma(f)$ instead of
$\Gamma(f,f)$.
An important property satisfied  by $\Gamma$ is that it is diffusive in the following sense:
\begin{equation}\label{diffu}
\Gamma(\phi(f),g)=\phi'(f)\Gamma(f,g).
\end{equation}
Besides, the eigenvalues of $-\mathcal{L}$
are given by
\[
\text{Sp}(-\mathcal{L})=\{i_1(i_1-1)+\cdots+i_{m}(i_{m}-1)\,|\,i_1,\ldots,i_m\in \nN\}.
\]
It may be ordered as a countable sequence like
$0 = \lambda_0 < \lambda_1 < \lambda_2 <  \cdots$, with a corresponding sequence of orthonormal eigenfunctions $u_0$, $u_1$, $u_2$, $\cdots$ where $u_0 =1$; in addition, this sequence of eigenfunctions forms a complete orthogonal basis of $L^2(\mu)$.
Also, note  that the first nonzero element of $\text{Sp}(-\mathcal{L})$
is $\lambda_1=1>0$.
Also, one can compute that, when $\lambda\in \text{Sp}(-\mathcal{L})$, then
$\text{Ker}(\mathcal{L}+\lambda\,I)$ is composed of those polynomial
functions $R(x_1,\ldots,x_m)$ having the form
\[
R(x_1,\ldots,x_m)=\sum_{i_1(i_1+1)+\cdots+i_{m}(i_{m}+1)=\lambda}\alpha(i_1,\cdots,i_{n_m}) J_{i_1}(x_1)\cdots J_{i_{m}}(x_{m}).
\]
Here $J_{i}(X)$ is the $i$th Jacobi polynomial, defined as
\[
J_{i}(x)=\frac{(-1)^i}{2^i i!}  \frac{d^i}{dx^i} \left\{(1 - x^2)^i\right\},\quad x\in\R.
\]
To end up with this quick summary, we recal the following Poincar\'e inequality, that is immediate to prove by using
the previous facts together with
the decomposition
$L^2(\mu) = \bigoplus_{\lambda\in \text{Sp}(-\mathcal{L})}
\text{Ker}(\mathcal{L}+\lambda\,I)$:
\begin{equation}\label{poincare}
{\rm Var}_\mu(f)\leq \int \Gamma(f)d\mu.
\end{equation}

{\bf Step 4}.
 We shall prove the existence of a constant $\kappa>0$,
depending on $p$, $\alpha$ and $d$ but \underline{not} on $n$, such
that, for any $\delta>0$,
\begin{equation}\label{claim3}
\sup_{n\geq 1} \,\,\,E_{\bf U}\left[ \frac{\delta}{\Gamma(Q_n)({\bf X})+\delta}
\right]{\bf 1}_{\{
{\rm Var}[Q_n({\bf X})|{\bf e},{\bf V}]\geq
\frac{d!}{2}\left(\frac{\alpha^{2}p}{3}\right)^d
\}}
\leq \kappa\,\delta^{\frac{1}{2d+1}}.
\end{equation}
The proof of (\ref{claim3}) will rely on the
Poincar\'e inequality (\ref{poincare}) which, here, takes the following form:
\begin{equation}\label{poincare2}
{\rm Var}[Q_n({\bf X})|{\bf e},{\bf V}]={\rm Var}_{\bf U}[Q_n({\bf X})]\leq E_{\bf U}[\Gamma(Q_n)({\bf X})].
\end{equation}
Another ingredient is the
Carbery-Wright inequality, that we recall for sake of completeness.

\begin{thm}(see \cite[Theorem 8]{CW})\label{cw-thm}
There exists an absolute constant $c>0$ such that, if $Q:\R^m\to\R$ is a polynomial of degree at most $k$
and $\mu$ is a log-concave probability measure on $\R^m$, then,
for all $\alpha>0$,
\begin{equation}\label{cw-ineq}
\left(\int Q^2d\mu\right)^{\frac1{2k}}\times
\mu\{x\in\R^m:\,|Q(x)|\leq \alpha\}
\leq c\,k\,\alpha^{\frac1k}.
\end{equation}
\end{thm}

Observe that the density of ${\bf U}$ is log-concave,
as an indicator function of a convex set.
Let us now proceed with the proof of (\ref{claim3}).
For any strictly positive $u$, and
provided ${\rm Var}[Q_n({\bf X})|{\bf e},{\bf V}]\geq
\frac{d!}{2}\left(\frac{\alpha^{2}p}{3}\right)^d$,
one has
\begin{equation}\label{cl45}
E\left[ \frac{\delta}{\Gamma(Q_n)({\bf X})+\delta}
\right]\leq \frac{\delta}{u} + P(\Gamma(Q_n)({\bf X})\leq u)
\leq \frac{\delta}{u} +c\, u^{\frac{1}{2d}},
\end{equation}
where $c>0$ denotes a constant only depending on $d$, $\alpha$ and $p$ and where the last inequality follows from the Carbery-Wright inequality (\ref{cw-ineq}), the inequality (\ref{poincare2}) and the fact that $\Gamma(Q_n)$
is a polynomial of order $2d$, see (\ref{gamma3}).
Finally, choosing $u=\delta^{\frac{2d}{2d+1}}$ in (\ref{cl45}) leads to the desired conclusion (\ref{claim3}).\\

{\bf Step 5}. We shall prove that
\begin{equation}\label{claim5}
\sup_{n\geq 1}
\left(
E_{\bf U} [\Gamma(\Gamma(Q_n))({\bf X})] + E_{\bf U}  \big|(\mathcal{L}Q_n)({\bf X})\big|
\right){\bf 1}_{\{
E[Q_n({\bf X})^2|{\bf e},{\bf V}]\leq M
\}}\leq c(M)
\end{equation}
where $c(M)$ is a constant only depending on $M$
(whose value may change from one line to another within this step).
First, relying on the results of Step 3 we have that, for any $n$,
\[
Q_n \in \bigoplus_{\alpha\leq \lambda_{2d}}\text{Ker}(\mathcal{L}+\alpha I).
\]
Since $\mathcal{L}$ is a bounded operator on the space $\bigoplus_{\alpha\leq \lambda_{2d}}\text{Ker}(\mathcal{L}+\alpha I)$, we deduce immediately that $\sup_{n\geq 1}
E_{\bf U} [(\mathcal{L}Q_n)({\bf X})^2]{\bf 1}_{\{E_{\bf U}[Q_n^2({\bf X})]\leq M\}}\leq c(M)$.
Besides, one has $\Gamma=\frac{1}{2}(\mathcal{L}+2\lambda I)$
on $\text{Ker}(\mathcal{L}+\lambda I)$ and one deduces for the same reason as above that
\[
\sup_{n\geq 1}E_{\bf U} [\Gamma(\Gamma(Q_n))({\bf X})]
{\bf 1}_{\{E_{\bf U}[Q_n^2({\bf X})]\leq M\}}\leq c(M).
\]
The proof of (\ref{claim5}) is complete.\\

{\bf Step 6}. We shall prove that, for any $n\geq 1$, any $0< \alpha\leq 1$, any $\delta>0$ and any $M>0$,
\begin{eqnarray}
\left|E_{\bf U}[(\phi-\phi\star p_\alpha)(Q_n({\bf X}))
]
\right|
{\bf 1}_{\{
{\rm Var}[Q_n({\bf X})|{\bf e},{\bf V}]\geq
\frac{d!}{2}\left(\frac{\alpha^{2}p}{3}\right)^d,
E[Q_n({\bf X})^2|{\bf e},{\bf V}]\leq M
\}}
\leq 2\kappa\,\delta^{\frac{1}{2d+1}}
+\sqrt{\frac{2}\pi}\,\frac{\alpha}{\delta}\,c(M).\notag\\
\label{claim4}
\end{eqnarray}
Using Step 4, one has
\begin{eqnarray}
&&\left|E_{\bf U}\left[
(\phi-\phi\star p_\alpha)(Q_n({\bf X}))
\frac{\Gamma(Q_{n})({\bf X})+\delta}{\Gamma(Q_{n})({\bf X})+\delta}
\right]
\right|
{\bf 1}_{\{
{\rm Var}[Q_n({\bf X})|{\bf e},{\bf V}]\geq
\frac{d!}{2}\left(\frac{\alpha^{2}p}{3}\right)^d,
E[Q_n({\bf X})^2|{\bf e},{\bf V}]\leq M
\}}
\notag\\
&\leq&
2\,E_{\bf U}\left[
\frac{\delta}{\Gamma(Q_{n})({\bf X})+\delta}
\right]{\bf 1}_{\{
{\rm Var}[Q_n({\bf X})|{\bf e},{\bf V}]\geq
\frac{d!}{2}\left(\frac{\alpha^{2}p}{3}\right)^d
\}}\notag\\
&&+
\left|E_{\bf U}\left[
(\phi-\phi\star p_\alpha)(Q_n({\bf X}))
\frac{\Gamma(Q_{n})({\bf X})}{\Gamma(Q_{n})({\bf X})+\delta}
\right]
\right|{\bf 1}_{\{
E[Q_n({\bf X})^2|{\bf e},{\bf V}]\leq M
\}}\notag
\\
&\leq&
2\,\kappa\,\delta^{\frac{1}{2d+1}}
+
\left|E_{\bf U}\left[
(\phi-\phi\star p_\alpha)(Q_n({\bf X}))
\frac{\Gamma(Q_{n})({\bf X})}{\Gamma(Q_{n})({\bf X})+\delta}
\right]
\right|{\bf 1}_{\{
E[Q_n({\bf X})^2|{\bf e},{\bf V}]\leq M
\}}.\label{FactC}
\end{eqnarray}
Now, set $\Psi(x)=\int_{-\infty}^{x} \phi(s) ds$ and let us apply (\ref{diffu}). We obtain
\begin{eqnarray}
\notag
&&\left|E_{\bf U}\left[
(\phi-\phi\star p_\alpha)(Q_n({\bf X}))
\frac{\Gamma(Q_{n})({\bf X})}{\Gamma(Q_{n})({\bf X})+\delta}
\right]
\right|\\
&=&\left|E_{\bf U}\left[\frac{1}{\Gamma(Q_{n})({\bf X})+\delta}\,\Gamma\big((\Psi-\Psi\star p_{\alpha})\circ Q_n,Q_n\big)({\bf X})\right]\right|\notag\\
\notag
&=&\left|E_{\bf U}\left[ ((\Psi-\Psi\star p_{\alpha})\circ Q_n)({\bf X}) \left(\Gamma\big(Q_n,\frac{1}{\Gamma(Q_{n})+\delta}\big)({\bf X})+\frac{(\mathcal{L}Q_n)({\bf X})}{\Gamma(Q_{n})({\bf X})+\delta}\right)
\right]\right|\\\notag
&=&\left|E_{\bf U}\left[((\Psi-\Psi\star p_{\alpha})\circ Q_n)({\bf X}) \left(-\frac{\Gamma(Q_n,\Gamma(Q_n))({\bf X})}
{(\Gamma(Q_{n})({\bf X})+\delta)^2}+\frac{(\mathcal{L}Q_n)({\bf X})}{\Gamma(Q_{n})({\bf X})+\delta}\right)\right]\right|\\
&\leq&\frac1\delta\,E_{\bf U}\left\{
|((\Psi-\Psi\star p_{\alpha})\circ Q_n)({\bf X})|\times\big[\Gamma(\Gamma(Q_n))({\bf X})+\big|(\mathcal{L}Q_n)({\bf X})\big|\big]\right\}.\label{FactA}
\end{eqnarray}
On the other hand, we have
\begin{eqnarray}
\notag
\left|\Psi(x)-\Psi\star p_{\alpha}(x)\right|&=&\left|\int_\R p_\alpha(y)\left(\int_{-\infty}^x \left(\phi(u)-\phi(u-y)\right)du\right)dy\right|\\\notag
&\leq&\int_\R p_\alpha(y)\left|\int_{-\infty}^x \phi(u)du-\int_{-\infty}^x\phi(u-y)du\right|dy\\
&\leq&\int_\R p_\alpha(y)\left|\int_{x-y}^x \phi(u)du\right| dy
\leq \int_\R p_{\alpha}(y)\left|y\right|dy \leq \sqrt{\frac{2}{\pi}}\alpha.\label{FactB}
\end{eqnarray}
The desired conclusion (\ref{claim4}) now follows easily from (\ref{claim5}), (\ref{FactC}), (\ref{FactA}) and (\ref{FactB}).\\

{\bf Step 7: Concluding the proof}.
Combining the results of all the previous steps,
we obtain, for any $n\geq 1$, any $0< \alpha\leq 1$, any $\delta>0$ and any $M>0$,
\begin{eqnarray}
&&\sup_{\phi\in C^\infty_c:\,\|\phi\|_\infty\leq 1}
\left|
E[\phi(Q_n({\bf X})]-E[\phi(Q_n({\bf Y})]
\right|\notag\\
&\leq&\frac{1}{\alpha}d_{FM}(Q_n({\bf X}),Q_n({\bf Y}))+\frac4M
+2P\left(
{\rm Var}[Q_n({\bf X})|{\bf e},{\bf V}]<
\frac{d!}{2}\left(\frac{\alpha^{2}p}{3}\right)^d
\right)\label{clbis}\\
&&+2P\left(
{\rm Var}[Q_n({\bf Y})|{\bf e},{\bf V}]<
\frac{d!}{2}\left(\frac{\alpha^{2}p}{3}\right)^d
\right)+4\kappa\,\delta^{\frac{1}{2d+1}}
+2\sqrt{\frac{2}\pi}\,\frac{\alpha}{\delta}\,c(M).
\label{cl}
\end{eqnarray}
In (\ref{clbis})-(\ref{cl}), take the limit $n\to\infty$.
Due to (\ref{cl-step1}) on one hand and Theorem \ref{MOO}
on the other hand (plus the fact that
the Fortet-Mourier distance $d_{FM}$ metrizes the convergence in distribution), one obtains
\[
\limsup_{n\to\infty}
\sup_{\phi\in C^\infty_c:\,\|\phi\|_\infty\leq 1}
\left|
E[\phi(Q_n({\bf X})]-E[\phi(Q_n({\bf Y})]
\right|\leq \frac4M +
4\kappa\,\delta^{\frac{1}{2d+1}}
+2\sqrt{\frac{2}\pi}\,\frac{\alpha}{\delta}\,c(M).
\]
The desired conclusion (\ref{cvdtv2})
then follows by letting (in this order) $\alpha\to 0$,
$\delta\to 0$ and $M\to\infty$.
\qed

\end{document}